\documentclass[10pt,leqno]{article}

\newcommand{\etalchar}[1]{$^{#1}$}
\def\polhk#1{\setbox0=\hbox{#1}{\ooalign{\hidewidth
  \lower1.5ex\hbox{`}\hidewidth\crcr\unhbox0}}}

\usepackage{amsmath,amstext,amssymb,amsfonts,amsthm}
\usepackage{ifthen,hyperref,dsfont}

\usepackage[a4paper,hscale=0.7,vscale=0.82,hmarginratio=1:1,includehead,headheight=14pt]{geometry}
\usepackage{fancyhdr}
\fancypagestyle{plain}{\fancyhead{}\fancyfoot{}}

\usepackage{color}

\newcommand{\wt}[1]{\ensuremath{\widetilde{#1}}}
\newcommand{\wh}[1]{\ensuremath{\widehat{#1}}}
\newcommand{\ol}[1]{\ensuremath{\overline{#1}}}
\newcommand{\fa}{\ensuremath{\text{\; for all \;}}}
\newcommand{\se}{\ensuremath{\subseteq}}

\newcommand{\cP}{\ensuremath{\mathcal{P}}}

\newcommand{\cH}{\ensuremath{\mathcal{H}}}

\newcommand{\R}{\ensuremath{\mathds{R}}}
\newcommand{\N}{\ensuremath{\mathds{N}}}
\newcommand{\Z}{\ensuremath{\mathds{Z}}}

\newcommand{\bH}{\ensuremath{\mathds{H}}}

\newcommand{\Diff}{\ensuremath{\operatorname{Diff}}}
\newcommand{\pr}{\ensuremath{\operatorname{pr}}}
\newcommand{\id}{\ensuremath{\operatorname{id}}}

\newcommand{\from}{\ensuremath{\nobreak:\nobreak}}
\renewcommand{\to}{\ensuremath{\nobreak\rightarrow\nobreak}}

\newcommand{\oU}[2]{U_{#1}^{[#2]}}
\newcommand{\cU}[2]{\ol{U}_{#1}^{[#2]}}
\newcommand{\sm}[1]{\wt{#1}}
\newcommand{\kq}[1]{\ensuremath{k_{#1}}}
\newcommand{\kp}[1]{\ensuremath{k_{#1}'}}
\newcommand{\supp}{\ensuremath{\operatorname{supp}}}

\theoremstyle{definition}
\newtheorem{definition}{Definition}[section]
\newtheorem{remark}[definition]{Remark}

\theoremstyle{plain}
\newtheorem{lemma}[definition]{Lemma}
\newtheorem{proposition}[definition]{Proposition}
\newtheorem{theorem}[definition]{Theorem}
\newtheorem{corollary}[definition]{Corollary}
\newtheorem*{nntheorem}{Theorem}

\newenvironment{prf}{\begin{proof}[\textbf{\upshape Proof.}]}{\end{proof}}

\newenvironment{introduction}
{\begin{center}{\textbf{\Large Introduction}}\end{center}
\markboth{Introduction}{Introduction}\noindent}{}

\begin{document}
\title{\textbf{Equivalences of Smooth and Continuous Principal Bundles
with In\-fi\-nite-\-Di\-men\-sion\-al Structure Group}}
\author{Christoph M\"uller%
\footnote{Funded in the DFG-Research Project 431/5-1,2:
\textit{Geometrische Darstellungstheorie wurzelgraduierter Lie-Gruppen}},
Christoph Wockel%
\footnote{Funded by a doctoral scholarship from the
\textit{Technische Universit\"at Darmstadt}}
\footnote{Present address: University of G\"ottingen, Mathematical Institute, Bunsenstr. 3-5, 37073 G\"ottingen, Germany}\\
        Fachbereich Mathematik\\
        Technische Universit\"at Darmstadt\\
\small\texttt{cmueller@mathematik.tu-darmstadt.de} \\
\small\texttt{christoph@wockel.eu}}
\maketitle
\thispagestyle{empty}

\begin{abstract} \noindent Let $K$ be a a Lie group, modeled on a
locally convex space, and $M$ a fi\-nite-\-di\-men\-sion\-al
paracompact manifold with corners. We show that each continuous
principal $K$-bundle over $M$ is continuously equivalent to a smooth
one and that two smooth principal $K$-bundles over $M$ which are
continuously equivalent are also smoothly equivalent. In the
concluding section, we relate our results to neighboring topics.
\\[\baselineskip] \noindent \textbf{Keywords:}
in\-fi\-nite-\-di\-men\-sion\-al Lie groups, manifolds with corners,
continuous principal bundles, smooth principal bundles, equivalences
of continuous and smooth principal bundles, smoothing of continuous
principal bundles, smoothing of continuous bundle equivalences, non
abelian \v{C}ech cohomology, twisted K-theory \\[\baselineskip]
\end{abstract}

\begin{introduction}%

This paper deals with the close interplay between continuous and
smooth principal $K$-bundles over $M$, where $K$ is a Lie group
modeled on an arbitrary locally convex space (following
\cite{milnor84}) and $M$ a fi\-nite-\-di\-men\-sion\-al paracompact
manifold with corners. In this paper we give a complete proof (of a
relative version) of the following theorem.

\begin{nntheorem}
Each continuous principal $K$-bundle over $M$ is equivalent to a
smooth principal $K$-bundle. Moreover, two smooth principal
$K$-bundles are continuously equivalent if and only if they are
smoothly equivalent.
\end{nntheorem}

One approach to a proof of this theorem is to introduce smooth
structures on classifying spaces and to smooth classifying maps. As an
example, the classifying space of $K=\operatorname{GL}_{n}$ is
isomorphic to the direct limit of the Grassmanians
\[
B\operatorname{GL}_{n}\cong \operatorname{Gr}_{n}(\infty) 
:=\lim_{\longrightarrow}\operatorname{Gr}_{n}(k).
\]
Then \cite[Th.\ 3.1]{gloeckner05} provides a smooth manifold structure
on $B\,\operatorname{GL}_{n}$, and one can smooth classifying maps as
in \cite[Th.\ 4.3.5]{hirsch76} for the case of vector bundles or in
Proposition
\ref{prop:smoothgingFiniteDimensionalBundlesViaClassifyingSpaces}, for
arbitrary fi\-nite-\-di\-men\-sion\-al principal bundles. In the
in\-fi\-nite-\-di\-men\-sion\-al case, the classifying space of the
diffeomorphism group $B\Diff(N)$ for a compact manifold $N$, which can
be viewed as a nonlinear Grassmanian, can also be given a smooth
structure \cite[44.21]{krieglmichor97}. 

Smooth structures on classifying spaces are considered in
\cite{mostow79}, but only generalized de Rham cohomologies are
constructed, and bundle theory is not discussed. However, a general
theory for differentiable structures on classifying spaces seems to be
missing. On the other hand, there exist partial answers to the above
question arising from the comparison of continuous and analytic fiber
bundles (cf.\ \cite{grauert58}, \cite{tognoli67} and
\cite{guaraldo02}). Since these considerations use strong constraints
on the structure group, e.g., its compactness in order to ensure a
smooth structure on its classifying space, they cannot be used in the
generality that we are aiming for.
\\

We now describe our results in some detail. In the first section, we
recall the basic facts on continuous and smooth principal bundles with
a focus on the description of bundles and bundle equivalences in terms
of locally trivial covers and transition functions. Furthermore, we
recall briefly the concept of differential calculus and the concept of
manifolds with corners that we use in this text.  In the end, we
outline how to prove our results for fi\-nite-\-di\-men\-sion\-al structure
groups by using smooth structures on classifying spaces.

The second section is exclusively devoted to the proofs of our main
results and to their technical prerequisites.  Lacking any smooth
structure on classifying spaces in general, we have to employ totally
different techniques coming from approximation results for Lie
group-valued functions (cf.\ Proposition
\ref{prop:smoothgingFiniteDimensionalBundlesViaClassifyingSpaces}). This
enables us to smooth representatives of continuous bundles or bundle
equivalences in combination with the fact that there is a large
freedom of choice in the description of principal bundles by locally
trivial covers and transition functions. In this way, we construct new
representatives of bundles and bundle equivalences that satisfy
cocycle or compatibility conditions on probably finer locally trivial
covers, but which describe equivalent objects. Since this technique
uses heavily the local compactness of the base manifold, there seems
to be no generalization of this method to
in\-fi\-nite-\-di\-men\-sion\-al base manifolds. Eventually, we
discover that the existence of smooth equivalent bundles and smooth
equivalences are a feature of convexity and continuity rendering
further requirements on $B\,K$ unnecessary.

In the third section, we relate our results to some neighboring
topics. In particular, we line out the relations to
\v{C}ech cohomology and to twisted $K$-theory. Concrete applications
of the above theorem arise, for instance, in twisted $K$-theory (cf.\
\cite{echterhoffEmersonKim06KKTheoreticDualityForProperTwistedActions},
\cite{mathaiStevenson06OnAGeneralisedCHKRTheorem}) and for obstruction
classes of lifting gerbes (cf.\ \cite{wagemannLaurent-Gengoux06}).

\end{introduction}

\section{Principal Fiber Bundles}

In this section we provide the basic material concerning manifolds
with corners and smooth and continuous principal bundles.

\begin{definition}[Continuous Principal Bundle]
\label{def:continuousBundle} Let $K$ be a topological group and $M$ be
a topological space. Then a \textit{continuous principal $K$-bundle
over $M$} (or shortly a \textit{continuous principal bundle}) is a
topological space $P$ together with a continuous right action $P\times K\to
P$, $(p,k)\mapsto p\cdot k$, and a map $\eta\from P\to M$ such that
there exists an open cover $(U_{i})_{i\in I}$ of $M$, called a \textit{locally
trivial cover}, and homeomorphisms
\[
\Omega_{i}:\eta^{-1}(U_{i})\to U_{i}\times K,
\]
called \textit{local trivializations}, satisfying $\pr_{1}\circ\;
\Omega_{i}=\left.\eta\right|_{\eta^{-1}(U_{i})}$ and $\Omega (p\cdot
k)=\Omega (p)\cdot k$. Here $K$ acts on $U_{i}\times K$ by right
multiplication in the second factor. We will use the calligraphic
letter $\cP$ for the tuple $(K,\eta\from P\to M)$.

A \textit{morphism of continuous bundles} or a \textit{continuous
bundle map} between two principal $K$-bundles $\cP$ and $\cP'$ over $M$ is
a continuous map $\Omega \from P\to P'$ satisfying $\Omega (p\cdot
k)=\Omega (p)\cdot k$.  Since $P'/K\cong M\cong P/K$, it induces a map
$\Omega^{\#}:M\to M$.  We call $\Omega$ a \textit{continuous bundle
equivalence} if it is an isomorphism and $\Omega^{\#}=\id_{M}$.
\end{definition}

\begin{remark}[Transition Functions]\label{rem:continuousBundle}
If $\cP$ is a continuous principal $K$-bundle over $M$, then the
local trivializations define continuous mappings $k_{ij}\from
U_{i}\cap U_{j}\to K$ by
\begin{equation}\label{eqn:transitionFunctions}
\Omega_{i}^{-1}(x,e)\cdot k_{ij}(x)=\Omega_{j}^{-1}(x,e)\fa x\in
U_{i}\cap U_{j},
\end{equation}
called \textit{transition functions}. The $k_{ij}$ satisfy the
\textit{cocycle condition}
\begin{equation}\label{eqn:cocycleCondition}
k_{ii}(x)=e\fa x\in U_{i}\quad \text{ and }\quad 
k_{ij}(x)\cdot k_{jn}(x)\cdot k_{ni}(x)=e\fa
x\in U_{i}\cap U_{j}\cap U_{n}.
\end{equation}
On the other hand, if $(V_{i})_{i\in I}$ is an open cover and
$k=(k_{ij})_{i,j\in I}$ is a collection of continuous maps
$k_{ij}:V_{i}\cap V_{j}\to K$ that satisfy condition
\eqref{eqn:cocycleCondition}, then
\[
P_{k}:=\bigcup_{j\in J}\{j\}\times V_{j}\times K /\sim \quad\text{ with
}\quad (j,x,k)\sim (j',x',k')\,\Leftrightarrow\, x=x' \text{ and }
k_{j'j}(x)\cdot k=k'\] defines a continuous principal $K$-bundle over
$M$. Here $\eta$ is given by $[i,x,k]\mapsto x$, the local
trivializations by $[(i,x,k)]\mapsto (x,k)$ and the $K$-action by
$([(i,x,k)],k')\mapsto [(i,x,kk')]$. We will write $\cP_{k}$ for a
bundle determined by a collection $(M,K,(V_{i})_{i\in
I},(k_{ij})_{i,j\in I})$.

If $k$ arises from the local trivializations of a given bundle
$\cP$ as in \eqref{eqn:transitionFunctions}, then
\[
\Omega \from P\to P_{k},\;\; p\mapsto [i,\Omega_{i}(p)]\text{ if }
p\in\eta^{-1}(U_{i})
\]
defines a bundle equivalence between $\cP$ and $\cP_{k}$ whose
inverse is given by $[i,x,k]\mapsto \Omega_{i}^{-1}(x,k)$.
\end{remark}

\begin{definition}[Differential Calculus on Locally Convex Spaces]
\label{def:differentialCalculus}(cf.\ \cite{gloecknerneeb})
Let $E$ and $F$ be locally convex spaces and $U\se E$ be open. Then
$f\from U\to F$ is called \textit{continuously differentiable} or
\textit{$C^{1}$} if it is continuous, for each $v\in E$ the
differential quotient
\[
df (x).v:=\lim_{h\to 0}\frac{1}{h}(f (x+hv)-f (x))
\]
exists and the map $df\from U\times E\to F$ is continuous. For $n>1$
we, recursively define 
\[
d^{n}f(x).(v_{1},\ldots,v_{n}):=\lim_{h\to 0}\frac{1}{h}
\left(d^{n-1}f(x+h).(v_{1},\ldots,v_{n-1})-
d^{n-1}f(x).(v_{1},\ldots,v_{n})\right)
\]
and say that $f$ is \textit{$C^{n}$} if $d^{k}f:U\times E^{k}\to F$
exists for all $k=1,\ldots,n$ and is continuous. We say that $f$ is
$C^{\infty}$ or \textit{smooth} if it is $C^{n}$ for all $n\in\N$.
\end{definition}

\begin{definition}[Lie Group]\label{def:lieGroup}
From the definition above, the notion of a \textit{Lie group} is
clear. It is a group which is a smooth manifold modeled on a locally
convex space such that the group operations are smooth.
\end{definition}

\begin{remark}[Convenient Calculus]
\label{rem:convenientCalculus}
We briefly recall the basic definitions
of the convenient calculus from \cite{krieglmichor97}. Again,
let $E$ and $F$ be locally convex spaces.  A curve $f:\R\to E$ is
called smooth if it is smooth in the sense of Definition
\ref{def:differentialCalculus}.  Then the $c^{\infty}$-topology on $E$
is the final topology induced from all smooth curves $f\in C^{\infty}
(\R,E)$. If $E$ is a Fr\'echet space, then the $c^{\infty}$-topology
is again a locally convex vector topology which coincides with the
original topology \cite[Th.\ 4.11]{krieglmichor97}.  If $U\se E$ is
$c^{\infty}$-open, then $f:U\to F$ is said to be of class $C^{\infty}$
or smooth if
\[
f_{*}\left(C^{\infty} (\R,U) \right)\se C^{\infty} (\R,F),
\]
i.e. if $f$ maps smooth curves to smooth curves. The chain rule
\cite[Prop.\ 1.15]{gloeckner02b} implies that each smooth map in the
sense of Definition \ref{def:differentialCalculus} is smooth in the
convenient sense. On the other hand, \cite[Th.\ 12.8]{krieglmichor97}
implies that on a Fr\'echet space a smooth map in the convenient sense
is smooth in the sense of Definition
\ref{def:differentialCalculus}. Hence for Fr\'echet spaces, the two
notions coincide.
\end{remark}

\begin{remark}[Manifold with Corners]\label{rem:manifoldWithCorners}
A \emph{$d$-dimensional manifold with corners} is a paracompact
topological space such that each point has a neighborhood that is
homeomorphic to an open subset of
\[
\R^{d}_{+}=\{(x_{1},\dotsc,x_{d})\in\R^{d}:x_{i}\geq 0\text{ for all
}i=1,\dotsc,d\}
\]
and such that the corresponding coordinate changes are smooth
(cf.\ \cite{smoothExt}, \cite{lee03}). The crucial point here is the notion of smoothness
for non-open domains. We define a map $f:A\se
\R^{n}\to \R^{m}$ to be smooth if for each $x\in A$, there exists a
neighborhood $U_{x}$ of $x$ which is open in $\R^{n}$, and a smooth
map $f_{x}:U_{x}\to \R^{m}$ such that $\left.f_{x}\right|_{A\cap
U_{x}}=\left.f\right|_{A\cap U_{x}}$.
\end{remark}

\begin{remark}[Paracompact Spaces]
\label{rem:topologicalPropertiesOfParacompactSpaces} We recall some
basic facts from general topology. If $X$ is a topological space, then
a collection of subsets $(U_{i})_{i\in I}$ of $X$ is called
\textit{locally finite} if each $x\in X$ has a neighborhood that has
non-empty intersection with only finitely many $U_{i}$, and $X$ is
called \textit{paracompact} if each open cover has a locally finite
refinement. If $X$ is the union of countably many compact subsets,
then it is called \emph{$\sigma$-compact}, and if each open cover has a
countable subcover, it is called \emph{Lindel\"of}.

Now let $M$ be a fi\-nite-\-di\-men\-sion\-al manifold with corners,
which is in particular locally compact and locally connected. For
these spaces, \cite[Theorems XI.7.2+3]{dugundji66} imply that $M$ is
paracompact if and only if each component is $\sigma$-compact,
equivalently, Lindel\"of.  Furthermore, \cite[Th.\
VIII.2.2]{dugundji66} implies that $M$ is normal in each of these
cases.
\end{remark}

\begin{definition}[Smooth Principal Bundle]\label{def:smoothBundle} If
$K$ is a Lie group and $M$ is a smooth manifold with corners, then a
continuous principal $K$-bundle over $M$ is called a \textit{smooth
principal $K$-bundle over $M$} (or shortly a \textit{smooth principal
bundle}) if the transition functions from Remark
\ref{rem:continuousBundle} are smooth for some choice of local
trivializations.
\end{definition}

\begin{remark}[Smooth Structure on Smooth Principal Bundles]
\label{rem:smoothStructureOnBundle} If $\cP$ is a smooth principal
bundle, then we define on $P$ the structure of a smooth manifold with
corners by requiring the local trivializations
\[
\Omega_{i}\from \eta^{-1}(U_{i})\to U_{i}\times K
\]
that define the smooth transition functions from Definition
\ref{def:smoothBundle} to be diffeomorphisms. This actually defines
a smooth structure on $P$, since it is covered by
$(\eta^{-1}(U_{i}))_{i\in I}$ and since the coordinate changes
\[
(U_{i}\cap U_{j})\times K\to (U_{i}\cap U_{j})\times K,\;\;
(x,k)\mapsto\Omega_{j}(\Omega_{i}^{-1}(x,k))=(x,k_{ij}(x)\cdot k)
\]
are smooth because the $k_{ij}$ are assumed to be smooth. A continuous
bundle map between smooth principal bundles is called a \textit{morphism of
smooth principal bundles} (or a \textit{smooth bundle map}) if it is
smooth with respect to the the smooth structure on the bundles just
described.
\end{remark}

\begin{remark}[Bundle Equivalences]\label{rem:bundleEquivalences}
If $\cP$ and $\cP'$ are two principal
$K$-bundles over $M$, then there exists an open cover $(U_{i})_{i\in
I}$ of $M$ such that we have local
trivializations
\begin{alignat*}{3}
\Omega_{i}&&\from \eta^{-1}(U_{i})&\to U_{i}\times K \\
\Omega'_{i}&&\from \eta'^{-1}(U_{i})&\to U_{i}\times K
\end{alignat*}
for $\cP$ and $\cP'$. In fact, if $(V_{j})_{j\in J}$ and
$(V'_{j'})_{j'\in J'}$ are locally trivial covers of $M$ (for $\cP$
and for $\cP'$, respectively), then 
\[
(V_{j}\cap V_{j'})_{(j,j')\in J\times J'}
\]
is simultaneously a locally trivial cover for both $\cP$ and $\cP '$,
and the local trivializations are given by restricting the
original ones.

If $\cP_{k}$ and $\cP_{k'}$ are given by transition
functions $k_{ij}$ and $k'_{ij}$ with respect to the same open cover
$(U_{i})_{i\in I}$ (i.e., $k_{ij}\from U_{i}\cap U_{j}\to K$ and $k'_{ij}\from
U_{i}\cap U_{j}\to K$), then a bundle equivalence $\Omega \from P_{k}\to
P_{k'}$ defines for each $i\in I$ a continuous map
\begin{equation}\label{eqn:definingTheBundleEquivalence}
\varphi_{i}\from U_{i}\times K\to K\text{\; by \;}
\Omega([(i,x,k)])=[(i,x,\varphi_{i}(x,k))].
\end{equation}
Furthermore, we have $\varphi_{i}(x,k)=\varphi_{i}(x,e)\cdot k$, since
$\Omega$ is assumed to satisfy $\Omega (p\cdot k)=\Omega (p)\cdot k$.
Setting $f_{i}(x):=\varphi_{i}(x,e)$, we thus obtain continuous maps
$f_{i}\from U_{i}\to K$ satisfying the \emph{compatibility condition}
\begin{equation}\label{eqn:compatibilityCondition}
f_{j}(x)=k'_{ji}(x)\cdot f_{i}(x)\cdot k_{ij}(x)\fa x\in
U_{i}\cap U_{j},
\end{equation}
since $[(i,x,k)]=[(j,x,k_{ji}(x)k)]$ has to be mapped to the same
element of $P_{k'}$ by $\Omega$. On the other hand, if for each
$i\in I$ we have continuous maps $f_{i}\from U_{i}\to K$ satisfying
\eqref{eqn:compatibilityCondition}, then
\[
P_{k}\ni [(i,x,k)]\mapsto [(i,x,f_{i}(x)\cdot k)]\in P_{k'}
\]
defines a bundle equivalence between $\cP_{k}$ and
$\cP_{k'}$ which covers the identity on $M$.

If $\cP_{k}$ and $\cP_{k'}$ are smooth and the maps $k_{ij}$
and $k_{ij}'$ are smooth, then it follows directly from
\eqref{eqn:definingTheBundleEquivalence} that a bundle equivalence
described by continuous maps $f_{i}:U_{i}\to K$ is smooth if and only
if these maps are smooth.
\end{remark}

\begin{lemma}[Smooth and Continuous Homotopies Coincide]
\label{lem:smoothingHomotopies}(\cite[Cor.\
12]{approx}, \cite{krieglMichor2002}) Let $M$ be a
fi\-nite-\-di\-men\-sion\-al manifold with corners and $N$ be a smooth
manifold modeled on a locally convex space. If $f:M\to N$ is
continuous, then there exists a continuous map $F:[0,1]\times M\to N$
such that $F(0,x)=f(x)$ and $F(1,\cdot)\from M\to N$ is
smooth. Furthermore, if $f,g:M\to N$ are smooth and there exists a
continuous homotopy between $f$ and $g$, then there exists a smooth
homotopy between $f$ and $g$.
\end{lemma}

\begin{lemma}[Smooth Structures on Classifying Spaces]
\label{lem:smoothStructureOnClassifyingSpacesForCompactLieGroups} If
$K$ is a compact Lie group, then it has a smooth classifying bundle
$EK\to BK$ (cf.\ \cite[Ch.\ 4.11]{husemoller94}), which is in general
infi\-nite-\-di\-men\-sion\-al.
\end{lemma}

\begin{prf}
Let $O_{k}\se \operatorname{GL}_{k}(\R)$ denote the orthogonal
group. If $k$ is sufficiently large, then we may identify $K$ with a
subgroup of $O_{k}$, and from \cite[Th.\ 19.6]{steenrod51} we get the
following formulae:
\begin{align*}
EK = & \lim_{\to} O_{n}/\big(O_{n-k}\times \id_{\R^{k}}\big),\\
BK = & \lim_{\to} O_{n}/\big(O_{n-k}\times K\big).
\end{align*}
Thus $EK$ and $BK$ are smooth manifolds by \cite[Th.\
3.1]{gloeckner05}, and since the action of $K$ is smooth, it follows
that $EK\to BK$ is a smooth $K$-principal bundle.
\end{prf}

\begin{proposition}[Smoothing Fi\-nite-\-Di\-men\-sion\-al Principal Bundles]
\label{prop:smoothgingFiniteDimensionalBundlesViaClassifyingSpaces} If
$\cP$ is a continuous principal $K$-bundle over $M$, $K$ is a
fi\-nite-\-di\-men\-sion\-al Lie group and $M$ is a fi\-nite-\-di\-men\-sion\-al manifold
with corners, then there exists a smooth bundle which is continuously
equivalent to $\cP$.  Moreover, two smooth principal $K$-bundles over
$M$ are smoothly equivalent if and only if they are continuously
equivalent.
\end{proposition}

\begin{prf}
Let $C$ be a maximal compact subgroup of $K$. Since $K/C$ is
contractible, there exists a $C$-reduction of $\cP$, i.e., we may choose
a locally trivial open cover $(U_{i})_{i\in I}$ with transition
functions $k_{ij}$ that take values in $C$.  They define a
continuous principal $C$-bundle which is given by a classifying map
$f\from M\to BC$.

By Lemma \ref{lem:smoothingHomotopies}, $f$ is homotopic to some
smooth map $\sm{f}\from M\to BC$ which in turn determines a smooth
principal $C$-bundle $\sm{\cP}$ over $M$ given by smooth transition
functions $\sm{k}_{ij}$. Furthermore, since $f$ and $\sm{f}$ are
homotopic, $\cP$ and $\sm{\cP}$ are equivalent, and we thus have a
continuous bundle equivalence given by continuous mappings
$f_{i}:U_{i}\to K$. The claim follows if we regard $k_{ij}$,
$\sm{k}_{ij}$ and $f_{i}$ as mappings into $K$.

Since smooth bundles yield smooth classifying maps and smooth
homotopies of classifying maps yield smooth bundle equivalences (all
the constructions in the topological setting depend only on partitions
of unity which we can assume to be smooth here), the second claim is
also immediate.
\end{prf}

\begin{remark}[On the previous proof]
The previous proof can also be obtained without the need of passing to
the direct limit in Proposition
\ref{lem:smoothStructureOnClassifyingSpacesForCompactLieGroups},
because $O_{n+k}/O_{n}$ is a ready a universal bundle if $\dim(M)\leq
n$ (cf.\ \cite[Rem.\ 19.7]{steenrod51}).
\end{remark}

\section{Equivalences of Smooth and Continuous Bundles}

In this section, we prove the two main results of this paper. We start
with the description of two important tools: a proposition for
smoothing continuous maps and a lemma for fading out continuous
functions. Then we provide some technical data for the proofs, namely
covers of the fi\-nite-\-di\-men\-sion\-al paracompact base manifold with
corners and suitable identity neighborhoods in the Lie group. On this
basis, we finally prove our claims after outlining the underlying
ideas in Remark \ref{rem:outlineOfTheProofs}.

\begin{remark}[Topology on $\boldsymbol{C(X,G)}$] If $X$ is a
Hausdorff space and $G$ is a topological group, then $C(X,G)_{c}$
denotes the topological group of continuous functions with respect to
pointwise multiplication and the topology of compact convergence. A
basic open identity neighborhood in this topology is given by
\[
\left\lfloor C,W \right\rfloor:=\{f\in C(X,G):f(C)\se W\}
\]
for a compact subset $C\se X$ and an open identity neighborhood $W\se G$.
\end{remark}

\begin{proposition}[Smoothing]\label{prop:smoothing}
Let $M$ be a fi\-nite-\-di\-men\-sion\-al manifold with corners, $K$ a Lie group
modeled on a locally convex space and $f\in C(M,K)$. If $A\subseteq M$
is closed and $U\subseteq M$ is open such that $f$ is smooth on a
neighborhood of $A\setminus U$, then each open neighborhood $O$ of $f$
in $C(M,K)_{c}$ contains a map $g$ which is smooth on a neighborhood
of $A$ and equals $f$ on $M\setminus U$.
\end{proposition}
\begin{prf}
This is \cite[Cor.\ 12]{approx}, also see \cite[Th.\
A.3.3]{neeb03} or \cite[Th.\ 2.5]{hirsch76}.
\end{prf}

\begin{remark}[Centered Chart, Convex Subset]
Let $K$ be a Lie group modeled on a locally convex topological vector
space $E$. A chart $\varphi:W\to\varphi(W)\subseteq E$ with $e\in W$
and $\varphi(e)=0$ is called a \emph{centered chart}. A subset $L$ of
$W$ is called \emph{$\varphi$-convex} if it is identified with a
convex subset $\varphi(L)$ in $E$. If $W$ itself is $\varphi$-convex,
we speak of a \emph{convex centered chart}.

It is clear that every open identity neighborhood in $K$ contains a
$\varphi$-convex open neighborhood for some centered chart $\varphi$,
because we can pull back any convex open neighborhood that is small
enough from the underlying locally convex vector space.
\end{remark}

\begin{lemma}[Fading-Out]\label{lem:fadingOut}
Let $M$ be a fi\-nite-\-di\-men\-sion\-al manifold with corners, $A$ and $B$ be
closed subsets satisfying $B\subseteq A^{0}$, $\varphi:W\to\varphi(W)$
be a convex centered chart of a Lie group $K$ modeled on a locally
convex space, and $f:A\to W$ be a continuous function. Then there is a
continuous function $F:M\to W\subseteq K$ that coincides with $f$ on
an open neighborhood of $B$ and is the identity on an open
neighborhood of $M\setminus A^{0}$. Moreover, $F$ can be chosen in a
way that if $W'\subseteq W$ is another $\varphi$-convex set containing
$e$, then $f(x)\in W'$ implies $F(x)\in W'$ for each $x\in A$, and if
$f$ is smooth on an open set $U\se A$, then $F$ is also smooth on $U$.
\end{lemma}

\begin{prf}
Since $M$ is paracompact and the closed sets $M\setminus A^{0}$ and
$B$ are disjoint by assumption, there exists a smooth map
$\lambda:M\to[0,1]$ such that $\lambda$ is identically $1$ on a
neighborhood of $B$ and is identically $0$ on a neighborhood of
$M\setminus A^{0}$ (see \cite[Th.  2.1]{hirsch76}). Since $\varphi(W)$
is a convex zero neighborhood in $E$, we have
$[0,1]\cdot\varphi(W)\subseteq\varphi(W)$. We use this to define the
continuous function
\begin{equation*}
f_{\lambda}:A\to W,\quad x\mapsto
\varphi^{-1}\Big(\lambda(x)\cdot\varphi\big(f(x)\big)\Big),
\end{equation*}
that coincides, by the choice of $\lambda$, with $f$ on $M\setminus
\supp(1-\lambda)\se B$ and is identically $e$ on
$M\setminus\supp(\lambda)\supseteq M\setminus A^{0}$. So we may extend
$f_{\lambda}$ to the continuous function
\begin{equation*}
F:M\to W,\quad x\mapsto\bigg\{\begin{array}{ll}
f_{\lambda}(x), & \text{if }x\in A \\
e, & \text{if }x\in M\setminus A^{0}
\end{array}
\end{equation*}
that satisfies all requirements.
\end{prf}

\begin{lemma}[Squeezing-in Manifolds with Corners]
\label{lem:sqeezingInManifoldsWithCorners} Let $W$ be an open
neighborhood of a point $x$ in $\R^{d}_{+}$ (cf.\ Remark
\ref{rem:manifoldWithCorners}) and $C\subseteq W$ be a compact set
containing $x$. Then there exists an open set $V$ satisfying $x\in
C\subseteq V\subseteq\ol{V}\subseteq W$ whose closure $\ol{V}$ is a
compact manifold with corners.
\end{lemma}
\begin{prf}
For every $x=(x_{1},\dotsc,x_{d})\in C$, there is an
$\varepsilon_{x}>0$ such that
\begin{equation}\label{eqn:cube}
B(x,\varepsilon):=[x_{1}-\varepsilon_{x},x_{1}+\varepsilon_{x}]\times
\dotsb\times[x_{d}-\varepsilon_{x},x_{d}+\varepsilon_{x}]\cap\R^{d}_{+}
\end{equation}
is contained in $W$.  The interiors $V_{x}:=B(x,\varepsilon_{x})^{0}$
in $\R^{d}_{+}$ form an open cover of the compact set $C$, of which we
may choose a finite subcollection $(V_{x_{i}})_{i=1,\dotsc,m}$
covering $C$. The union $V:=\bigcup_{i=1}^{m}V_{x_{i}}$ satisfies all
requirements. In particular, $\ol{V}$ is a compact manifold with
corners, because it is a finite union of cubes whose sides are
orthogonal to the coordinate axes.
\end{prf}

\begin{proposition}[Nested Covers]\label{prop:NestedCovers}
Let $M$ be a connected fi\-nite-\-di\-men\-sion\-al manifold with corners and
$(U_{j})_{j\in J}$ be an open cover of $M$. Then there exist countable
open covers $\big(\oU{i}{\infty}\big)_{i\in\N}$ and
$\big(\oU{i}{0}\big)_{i\in\N}$ of $M$ such that
$\cU{i}{\infty}:=\ol{\oU{i}{\infty}}$ and $\cU{i}{0}:=\ol{\oU{i}{0}}$
are compact manifolds with corners, $\cU{i}{\infty}\subseteq\oU{i}{0}$
for all $i\in\N$, and such that even the cover
$\big(\cU{i}{0}\big)_{i\in\N}$ of $M$ by compact sets is locally
finite and subordinate to $(U_{j})_{j\in J}$.

In this situation, let $L$ be any countable subset of the open
interval $(0,\infty)$. Then for every $\lambda\in L$, there exists a
countable, locally finite cover $\big(\oU{i}{\lambda}\big)_{i\in\N}$
of $M$ by open sets whose closures are compact manifolds with corners
such that $\cU{i}{\lambda}\subseteq\oU{i}{\mu}$ holds whenever
$0\leq\mu<\lambda\leq\infty$.
\end{proposition}
\begin{prf}
For every $x\in M$, we have $x\in U_{j(x)}$ for some $j(x)\in J$. Let
$(U_{x},\varphi_{x})$ be a chart of $M$ around $x$ such that
$\ol{U_{x}}\subseteq U_{j(x)}$. We can even find an open neighborhood
$V_{x}$ of $x$ whose closure $\ol{V_{x}}$ is compact and contained in
$U_{x}$. Since $M$ is paracompact, the open cover
$\big(V_{x}\big)_{x\in M}$ has a locally finite subordinated cover
$(V_{i})_{i\in I}$, where $V_{i}\subseteq V_{x}$ and
$\ol{V_{i}}\subseteq\ol{V_{x}}\subseteq U_{x}$ for suitable
$x=x(i)$. Since $M$ is also Lindel\"of, we may assume that $I=\N$.

To find suitable covers $\oU{i}{\infty}$ and $\oU{i}{0}$, we are going
to enlarge the sets $V_{i}$ so carefully in two steps that the
resulting covers remain locally finite. More precisely,
$\oU{i}{\infty}$ and $\oU{i}{0}$ will be defined inductively so that
even the family $(V_{k}^{i})_{k\in\N}$ with
\begin{equation*}
V_{k}^{i}:=\left\{\begin{array}{ll}
\cU{k}{0} & \text{ for }k\leq i \\
V_{k} & \text{ for }k>i
\end{array}\right.
\end{equation*}
is still a locally finite cover of $M$ for every $i\in\N_{0}$. We
already know this for $i=0$, because $V_{k}^{0}=V_{k}$ for all
$k\in\N$. For $i>0$, we proceed by induction.

For every point $y\in\ol{V_{i}}$, there is an open neighborhood
$V_{i,y}$ of $y$ inside $U_{x(i)}$ whose intersection with just
finitely many $V^{i-1}_{j}$ is non-empty. Under the chart
$\varphi_{x(i)}$, this neighborhood $V_{i,y}$ is mapped to an open
neighborhood of $\varphi_{x(i)}(y)$ in the modeling space $\R^{d}_{+}$
of $M$. There exist real numbers
$\varepsilon_{0}(y)>\varepsilon_{\infty}(y)>0$ such that the cubes
$B(y,\varepsilon_{\infty}(y))$ and $B(y,\varepsilon_{0}(y))$
introduced in \eqref{eqn:cube} are compact neighborhoods of
$\varphi_{x(i)}(y)$ contained in $\varphi_{x(i)}(V_{i,y})$. Since
$\ol{V_{i}}$ is compact, it is covered by finitely many sets
$V_{i,y}$, say by $\big(V_{i,y})_{y\in Y}$ for a finite subset $Y$ of
$\ol{V_{i}}$.  We define the open sets
\begin{equation*}
\oU{i}{\infty}:=\bigcup_{y\in Y}
\varphi_{x(i)}^{-1}\left(B(y,\varepsilon_{\infty}(y))^{0}\right)
\text{ and }
\oU{i}{0}:=\bigcup_{y\in Y}
\varphi_{x(i)}^{-1}\left(B(y,\varepsilon_{0}(y))^{0}\right),
\end{equation*}
whose closures are compact manifolds with corners, because they are a
finite union of cubes under the chart $\varphi_{x(i)}$. On the one
hand, the construction guarantees
\begin{equation*}
V_{i}\subseteq\oU{i}{\infty}\subseteq\cU{i}{\infty}\subseteq
\oU{i}{0}\subseteq\cU{i}{0}\subseteq\bigcup_{y\in Y}V_{i,y}
\subseteq U_{x(i)}.
\end{equation*}
On the other hand, the cover $\big(V^{i}_{k}\big)_{k\in\N}$ is locally
finite, because it differs from the locally finite cover
$\big(V^{i-1}_{k}\big)_{k\in\N}$ in the single set
$V_{i}^{i}=\cU{i}{0}$.

For a proof of the second claim, we fix an enumeration
$\lambda_{1},\lambda_{2},\dotsc$ of $L$ for an inductive construction
of the covers.  Then for any $n\geq 1$ and $i\in\N$, we apply Lemma
\ref{lem:sqeezingInManifoldsWithCorners} to
$C:=\varphi_{i}\big(\cU{i}{\ol{\lambda}}\big)$ and
$W:=\varphi_{i}\big(\oU{i}{\underline{\lambda}}\big)$, where
$\ol{\lambda}$ (resp. $\underline{\lambda}$) is the smallest
(resp. largest) element of $\lambda_{1},\dotsc,\lambda_{n-1}$ larger
than (resp. smaller than) $\lambda_{n}$ for $n>1$ and $\infty$
(resp. $0$) for $n=1$. We get open sets $\oU{i}{\lambda_{n}}$ such
that the condition $\cU{i}{\lambda}\subseteq\oU{i}{\mu}$ holds
whenever $0\leq\mu<\lambda\leq\infty$ are elements in
$\{\lambda_{1},\dotsc,\lambda_{n}\}$, and eventually in $L$. This
completes the proof.
\end{prf}

\begin{remark}[Locally Finite Covers by Compact Sets]
\label{rem:locallyFiniteCoversByCompactSets} If
$\big(\ol{U}_{i}\big)_{i\in I}$ is a locally finite cover of $M$ by
compact sets, then for fixed $i\in I$, the intersection
$\ol{U}_{i}\cap\ol{U}_{j}$ is non-empty for only finitely many $j\in
I$. Indeed, for every $x\in\ol{U}_{i}$, there is an open neighborhood
$U_{x}$ of $x$ such that $I_{x}:=\{j\in
I:U_{x}\cap\ol{U}_{j}\not=\emptyset\}$ is finite. Since $\ol{U}_{i}$
is compact, it is covered by finitely many of these sets, say by
$U_{x_{1}},\dotsc,U_{x_{n}}$. Then $J:=I_{x_{1}}\cup\dotsb\cup
I_{x_{n}}$ is the finite set of indices $j\in J$ such that
$\ol{U}_{i}\cap\ol{U}_{j}$ is non-empty, proving the claim.
\end{remark}

\begin{remark}[Intersections]
From now on, multiple lower indices on subsets of $M$ always indicate
intersections, namely $U_{1\dotsb r}:=U_{1}\cap\dotsc\cap U_{r}$.
\end{remark}

\begin{lemma}[Suitable Identity Neighborhoods]
\label{lem:SuitableIdentityNeighborhoods} Let $M$ be a
fi\-nite-\-di\-men\-sion\-al manifold with corners that is covered locally
finitely by countably many compact sets
$\big(\ol{U}_{i}\big)_{i\in\N}$. Moreover, let $k_{ij}:\ol{U}_{ij}\to
K$ be continuous functions into a Lie group $K$ so that
$k_{ij}=k_{ji}^{-1}$ holds for all $i,j\in\N$. Then for any convex
centered chart $\varphi:W\to\varphi(W)$ of $K$, each sequence of open
unit neighborhoods $(W'_{j})_{j\in \N}$ with $W'_{j}\se W$ and each
$\alpha \in \N $, there are $\varphi$-convex open identity
neighborhoods $W_{ij}^{\alpha}\subseteq W$ in $K$ for indices $i<j$
and $W_{j}^{\alpha}\se W'_{j}$ for $j\in\N$ that satisfy
\begin{align}
k_{ji}(x)\cdot(W_{ij}^{\alpha})^{-1}\cdot W_{i}^{\alpha}\cdot 
k_{ij}(x)\subseteq W_{j}^{\alpha}
&\fa x\in\ol{U}_{ij\alpha}\text{ and } i<j,
\label{eqn:newf} \\
k_{ji}(x)\cdot(W_{ij}^{\alpha})^{-1}\cdot W_{in}^{\alpha}\cdot 
k_{ij}(x)\subseteq W_{jn}^{\alpha}
&\fa x\in\ol{U}_{ijn\alpha}\text{ and } i<j<n
\label{eqn:newrec}
\end{align}
\end{lemma}
\begin{prf}
Inititally, we set $W_{i}^{\alpha}:=W'_{i}$ for all $i$, respectively
$W_{ij}^{\alpha}:=W$ for all $i<j$, disregarding the conditions
\eqref{eqn:newf} and \eqref{eqn:newrec}. These sets are
shrinked later so that they satisfy \eqref{eqn:newf} and
\eqref{eqn:newrec}.

Our first goal is to satisfy \eqref{eqn:newf}. We note that the
condition in \eqref{eqn:newf} becomes trivial if $\ol{U}_{j\alpha}$ is
empty, because this implies $\ol{U}_{ij\alpha}=\emptyset$. So we need
to consider at most finitely many conditions \eqref{eqn:newf}
corresponding to the finitely many $j\in\N$ such that
$\ol{U}_{j\alpha}\neq\emptyset$, and we deal with those
inductively in decreasing order of $j$, starting with the maximal
such index.

For fixed $j$ and all $i<j$ with $\ol{U}_{ij\alpha}\neq\emptyset$, we
describe below how to make the $\varphi$-convex identity neighborhoods
$W_{ij}^{\alpha}$ and $W_{i}^{\alpha}$ on the left hand side smaller
so that the corresponding conditions \eqref{eqn:newf} are satisfied.
Making $W_{ij}^{\alpha}$ and $W_{i}^{\alpha}$ smaller does not
compromise any conditions on $W_{ij'}^{\alpha}$ and $W_{j'}^{\alpha}$
for $j'>j$ that we guaranteed before, because these sets can only
appear on the left hand side of such conditions.

To satisfy condition \eqref{eqn:newf} for given $i<j$ and
$W_{j}^{\alpha}$, we note that the function
\begin{equation*}
\varphi_{ij}\from 
\ol{U}_{ij\alpha}\times K\times K\to K,
\quad(x,k,k')\mapsto k_{ji}(x)\cdot k^{-1}\cdot k'\cdot k_{ij}(x)
\end{equation*}
is continuous and maps all the points $(x,e,e)$ for
$x\in\ol{U}_{ij\alpha}$ to the identity $e$ in $K$. Hence we may
choose open neighborhoods $U_{x}$ of $x$ and $\varphi$-convex open
identity neighborhoods $W_{x}\subseteq W_{ij}^{\alpha}$ and
$W_{x}'\subseteq W_{i}^{\alpha}$ such that $\varphi_{ij}(U_{x}\times
W_{x}\times W_{x}')\subseteq W_{j}^{\alpha}$. Since
$\ol{U}_{ij\alpha}$ is compact, it is covered by finitely many
$U_{x}$, say by $(U_{x})_{x\in F}$ for a finite set $F\se
\ol{U}_{ij\alpha}$. Then we replace $W_{ij}^{\alpha}$ and
$W_{i}^{\alpha}$ by their subsets $\bigcap_{x\in F}W_{x}$ and
$\bigcap_{x\in F}W_{x}'$, respectively, which are $\varphi$-convex
open identity neighborhoods such that
$\varphi_{ij}(\ol{U}_{ij\alpha}\times W_{ij}^{\alpha}\times
W_{i}^{\alpha})\subseteq W_{j}^{\alpha}$, in other words,
\eqref{eqn:newf} is satisfied

Our second goal is to make the sets $W_{ij}^{\alpha}$ also satisfy
\eqref{eqn:newrec}, which is non-trivial for the finitely many triples
$(i,j,n)\in\N^{3}$ with $i<j<n$ that satisfy
$\ol{U}_{ijn\alpha}\neq\emptyset$. We can argue as above,
except for a slightly more complicated order of processing the sets
$W_{jn}^{\alpha}$ on the right hand side. Namely, we define the
following total order
\begin{equation}\label{eqn:order}
(i,j)<(i',j')\quad:\Leftrightarrow\quad
j<j'\text{ or }(j=j'\text{ and }i<i')
\end{equation}
on pairs of real numbers, in particular on pairs of indices $(i,j)$ in
$\N\times\N$ with $i<j$. Note that this guarantees $(i,j)<(j,n)$ and
$(i,n)<(j,n)$ whenever $i,j,n$ are as in condition
\eqref{eqn:newrec}. We process the pairs $(j,n)$ with
$\ol{U}_{ijn\alpha}\neq\emptyset$ for some $i$ in descending order,
starting with the maximal such pair. At each step, we fix
$W_{jn}^{\alpha}$ and make $W_{ij}^{\alpha}$ and $W_{in}^{\alpha}$
smaller for all relevant $i<j$ so that \eqref{eqn:newrec} is
satisfied. This does not violate any conditions \eqref{eqn:newf} or
\eqref{eqn:newrec} that we guaranteed earlier in the process, because
$W_{ij}^{\alpha}$ and $W_{in}^{\alpha}$ can only appear on the left
hand side of such conditions. For the choice of the smaller identity
neighborhoods, we use the continuous function
\begin{equation*}
\varphi_{ijn}\from\ol{U}_{ijn\alpha}\times K\times
K \to K,\quad
(x,k,k')\mapsto k_{ji}(x)\cdot k^{-1}\cdot k'\cdot k_{ij}(x)
\end{equation*}
and the compactness of $\ol{U}_{ijn\alpha}$ and argue as before. We
thus accomplish our second goal.
\end{prf}

\begin{remark}[Outline of the Proofs]\label{rem:outlineOfTheProofs}
Although the proofs of our main results are quite technical, the
underlying ideas are easy to explain. The following two theorems
require us to construct principal bundles and/or equivalences between
them, and we always do so locally on countable covers of the base
manifold by induction. In these constructions, every new transition
function (respectively, every new local representative of an
equivalence)
\begin{itemize}
\item is already determined by cocycle conditions (respectively, by
compatibility conditions) on a subset of its domain,
\item from which it will be ``faded out'' to a continuous function on
the whole domain
\item and smoothed, if necessary.
\end{itemize}
In each such step, we need a safety margin to modify the functions
without compromising previous achievements too much, and these safety
margins are the nested open covers provided by Proposition
\ref{prop:NestedCovers}. In order to ``fade out'' appropriately, we
need to make sure that the values of the corresponding functions stay
in certain identity neighborhoods of the structure group. This is
achieved with the data from Lemma
\ref{lem:SuitableIdentityNeighborhoods}.
\end{remark}

\begin{theorem}[Smoothing Continuous Principal Bundles]
\label{thm:smoothingContinuousPrincipalBundles} Let $K$ be a Lie group
modeled on a locally convex space, $M$ be a
fi\-nite-\-di\-men\-sion\-al connected paracompact manifold with
corners and $\cP$ be a continuous principal $K$-bundle over $M$. If
$C\se M$ is closed and the restriction of $\cP$ to some open
neighborhood of $C$ is smooth, then there exists an open neighborhood
$T$ of $C$ such that the restriction $\left.\cP\right|_{T}$ extends to
a smooth principal $K$-bundle $\sm{\cP}$ over $M$, in the sense that
$\left.\cP\right|_{T}$ is a $K$-invariant open subset of $\sm{\cP}$
and the $K$-action and bundle projection from $\left.\cP\right|_{T}$
extend to the ones on $\sm{\cP}$. Furthermore, there exists a
continuous bundle equivalence $\Omega \from P\to\sm{P}$, which
restricts to the identity on $\eta^{-1}(T)$.
\end{theorem}

\begin{prf}
We assume that the continuous bundle $\cP$ is given by $\cP_{k}$ as in
Remark \ref{rem:continuousBundle}, where $(U_{j})_{j\in J}$ is a
locally trivial cover of $M$ and $k_{ij}:U_{ij}\to K$ are continuous
transition functions that satisfy the cocycle condition $k_{ij}\cdot
k_{jn}=k_{in}$ pointwise on $U_{ijn}$. That $\cP$ is smooth on a
neighborhood of $C$ implies that there exists an open neighborhood $S$
of $C$ such that the restriction of each $k_{ij}$ to $U_{ij}\cap S$ is
smooth. In fact, let $S'$ be an open neighborhood of $C$ such that
$\left.\cP\right|_{S'}$ is smooth. Since $M$ is normal (see Remark
\ref{rem:topologicalPropertiesOfParacompactSpaces}), we find open sets
$S,T\se M$ that satisfy $C\se T\se \ol{T}\se S\se \ol{S}\se S'$, which
we fix from now on. In addition, there exists a locally trivial cover
of $S'$, together with local trivialisations, such that the resulting
transition functions are smooth. Restricting the continuous transition
functions of an arbitrary locally trivial cover to the complement of
$\ol{S}$, adding the smooth ones and the ones induced by the cocycle
condition yields the desired collection of transition functions.

Proposition \ref{prop:NestedCovers} yields open covers
$\big(\oU{i}{\infty}\big)_{i\in\N}$ and $\big(\oU{i}{0}\big)_{i\in\N}$
of $M$ subordinate to $(U_{j})_{j\in J}$ with
$\cU{i}{\infty}\se\oU{i}{0}$ for all $i\in\N$. For every $i\in\N$, we
denote by $U_{i}$ an open set of the cover $(U_{j})_{j\in J}$ that
contains $\cU{i}{0}$ and observe that $(U_{i})_{i\in\N}$ is still a
locally trivial open cover of $M$. In our construction, we need open
covers not only for pairs $(j,n)\in\N\times\N$ with $j<n$, but also
for pairs $(j-1/3,n)$, $(j-2/3,n)$ in-between and $(n,n)$ to enable
continuous extensions and smoothing. The function
\begin{equation*}
\lambda:\left\{(j,n)\in\frac{1}{3}\N_{0}\times\N:j\leq n\right\}
\to[0,\infty),\quad\lambda(j,n)=\frac{n(n-1)}{2}+j,
\end{equation*}
is tailored to map the pairs
$(0,1),(1,1),(1,2),(2,2),(1,3),(2,3),(3,3),(1,4),\dotsc$ to the
integers $0,1,2,\dotsc$, respectively, and the other pairs
in-between. If we apply the second part of Proposition
\ref{prop:NestedCovers} to the countable subset
$L:=(\operatorname{im}\lambda)\setminus\{0\}$ of $(0,\infty)$, we get
open sets $\oU{i}{jn}:=\oU{i}{\lambda(j,n)}$ for all pairs $(j,n)$ in
the domain of $\lambda$ such that $\big(\cU{i}{jn}\big)_{i\in\N}$ are
again locally finite covers. We note that $(j,n)<(j',n')$ in the sense
of \eqref{eqn:order} implies $\cU{i}{j'n'}\subseteq\oU{i}{jn}$.

Let $\varphi:W\to\varphi(W)$ be an arbitrary convex centered chart of
$K$ and consider the countable compact cover
$\big(\cU{i}{0}\big)_{i\in\N}$ of $M$ and the restrictions
$k_{ij}|_{\cU{ij}{0}}$ of the continuous transition functions to the
corresponding intersections. Then Lemma
\ref{lem:SuitableIdentityNeighborhoods}, applied to the sequence of
open unit neighborhoods which is constantly $W$, yields open
$\varphi$-convex identity neighborhoods $W_{ij}^{\alpha}$ and
$W_{i}^{\alpha}$ with the corresponding properties.

Our first goal is the construction of smooth maps
$\sm{k}_{ij}:\oU{ij}{0}\to K$ that satisfy the cocycle condition on
the open cover $\big(\oU{i}{\infty}\big)_{i\in\N}$ of $M$, which
uniquely determines a smooth principal $K$-bundle $\cP_{\sm{k}}$
by Remarks \ref{rem:continuousBundle} and
\ref{rem:smoothStructureOnBundle}. Furthermore, we shall construct
$\sm{k}_{ij}$ in a way that guarantees
\[
\sm{k}_{ij}\Big|_{\oU{ij}{0}\cap T}=k_{ij}\Big|_{\oU{ij}{0}\cap T}
\]
ensuring that we may view $\left.\cP_{\sm{k}}\right|_{T}$ as a subset
of $\cP_{k}$.  These maps $\sm{k}_{ij}$ will be constructed
step-by-step in increasing order with respect to \eqref{eqn:order},
starting with the minimal index $(1,2)$. At all times during the
construction, the conditions
\begin{enumerate}
\renewcommand{\labelenumi}{(\theenumi)}
\renewcommand{\theenumi}{\alph{enumi}}
\item $\sm{k}_{jn}=\sm{k}_{ji}\cdot\sm{k}_{in}$ pointwise on
$\cU{ijn}{jn}$ for all $i<j<n$ in $\N$,\label{eqn:cocycle}
\item $\big(\sm{k}_{jn}\cdot k_{nj}\big)\big(\cU{jn\alpha}{jn}\big)
\subseteq W_{jn}^{\alpha}$ for all $j<n$ and $\alpha$ in
$\N$ and \label{eqn:quotient}
\item $\sm{k}_{jn}\Big|_{\oU{jn}{0}\cap T}=k_{jn}\Big|_{\oU{jn}{0}\cap T}$
for all $j,n\in\N$\label{eqn:relative1}
\end{enumerate}
will be satisfied whenever all $\sm{k}_{ij}$ involved have already
been constructed. We are now going to construct the smooth maps
$\sm{k}_{jn}$ for indices $j<n$ in $\N$ (and implicitly $\sm{k}_{nj}$
as $\sm{k}_{nj}(x):=\sm{k}_{jn}(x)^{-1}$), assuming that this has
already been done for pairs of indices smaller than $(j,n)$.
\begin{itemize}

\item To satisfy all relevant cocycle conditions, we start with
\begin{equation*}
\sm{k}'_{jn}:\bigcup_{i<j}\cU{ijn}{j-1,n}\to K,\quad
\sm{k}'_{jn}(x):=\sm{k}_{ji}(x)\cdot\sm{k}_{in}(x)\text{ for }
x\in\cU{ijn}{j-1,n}.
\end{equation*}
This function is well-defined, because the cocycle conditions
\eqref{eqn:cocycle} for lower indices assert that for any indices
$i'<i<j$ and any point $x\in\cU{i'jn}{j-1,n}\cap\cU{ijn}{j-1,n}$, we
have
\begin{equation*}
\sm{k}_{ji'}(x)\cdot\sm{k}_{i'n}(x)
=\sm{k}_{ji'}(x)\cdot\sm{k}_{i'i}(x)
\cdot\sm{k}_{ii'}(x)\cdot\sm{k}_{in}(x)
=\sm{k}_{ji}(x)\cdot\sm{k}_{in}(x),
\end{equation*}
because $\cU{i'ijn}{j-1,n}$ is contained in both $\cU{i'ij}{ij}$ and
$\cU{i'in}{in}$. Furthermore, $\sm{k}'_{jn}$ coincides with $k_{jn}$
on $\bigcup_{i<j}\cU{ijn}{j-1,n}\cap T$ as
$\sm{k}'_{jn}(x)=\sm{k}_{ji}(x)\cdot \sm{k}_{in}(x)=k_{ji}(x)\cdot
k_{in}(x)=k_{jn}(x)$.

\item Next, we want to extend the map $\sm{k}'_{jn}$ on
$\bigcup_{i<j}\cU{ijn}{j-1,n}$ to a continuous map $k'_{jn}$ on
$\oU{jn}{0}$ without compromising the cocycle conditions too much. To
do this, we consider the function
$\varphi_{jn}:=\sm{k}_{jn}'\cdot k_{nj}:\bigcup_{i<j}^{}\cU{ijn}{j-1,n}\to
K$. For all $i<j$, $\alpha\in\N$ and $x\in\cU{ijn\alpha}{j-1,n}$,
conditions \eqref{eqn:quotient} above and \eqref{eqn:newrec} of Lemma
\ref{lem:SuitableIdentityNeighborhoods} imply
\begin{align*}
\varphi_{jn}(x)
&=(\sm{k}'_{jn}k_{nj})(x)
=k_{ji}(x)\cdot
\big(\underbrace{(\sm{k}_{ij}\cdot k_{ji})(x)}_{\in W_{ij}^{\alpha}}\big)^{-1}
\cdot\underbrace{(\sm{k}_{in}\cdot k_{ni})(x)}_{\in W_{in}^{\alpha}}
\cdot k_{ij}(x) \\
&\in k_{ji}(x)\cdot(W_{ij}^{\alpha})^{-1}\cdot W_{in}^{\alpha}
\cdot k_{ij}(x)
\subseteq W_{jn}^{\alpha},
\end{align*}
because $\cU{ijn\alpha}{j-1,n}$ is contained in both
$\cU{ij\alpha}{ij}$ and $\cU{in\alpha}{in}$. Since the values of
$\varphi_{jn}$ are contained in particular in the identity
neighborhood $W$, we may apply Lemma \ref{lem:fadingOut} to
$M:=\oU{jn}{0}$ and its subsets $A:=\bigcup_{i<j}\cU{ijn}{j-1,n}$ and
$B:=\bigcup_{i<j}\cU{ijn}{j-2/3,n}$. This yields a continuous function
$\Phi_{jn}\from \oU{jn}{0}\to W$ that coincides with $\varphi_{jn}$ on
a neighborhood of $B$, is the identity on a neighborhood of
$\oU{jn}{0}\setminus A^{0}$ and satisfies $\Phi_{jn}(x)\in
W_{jn}^{\alpha}$ for all $x\in\cU{jn\alpha}{j-1,n}$. We define
$k'_{jn}\from \oU{jn}{0}\to K$ by $k'_{jn}:=\Phi_{jn}k_{jn}$ and note
that $k'_{jn}$ coincides with the smooth function $\sm{k}'_{jn}$ on a
neighborhood of $B$ and with $k_{jn}$ on a neighborhood of
$\oU{jn}{0}\setminus A^{0}$.  In addition, $\Phi_{jn}$ is the identity on
$A\cap T$, because $\sm{k}'_{jn}$ and $k_{jn}$ coincide
there. Furthermore, by the last conclusion from Lemma
\ref{lem:fadingOut}, $\Phi_{jn}$ is smooth on $\oU{jn}{0}\cap S$, because
$\sm{k}'_{jn}$ and $k_{nj}$ are smooth on $A^{0}\cap S$ and
$\Phi_{jn}$ is the identity on an open neighborhood of
$\oU{jn}{0}\setminus A^{0}$. Consequently, $k'_{jn}$ coincides with
$k_{jn}$ on $\oU{jn}{0}\cap T$ and is smooth on $\oU{jn}{0}\cap S$,
which is an open neighborhood of $\oU{jn}{0}\cap \ol{T}$.

\item We finally get the smooth map $\sm{k}_{jn}\from \oU{jn}{0}\to K$
that we are looking for if we apply Proposition \ref{prop:smoothing}
to the function $k'_{jn}$ on $M:=A:= \oU{jn}{0}$, to the open
complement $U$ of $\bigcup_{i<j}\cU{ijn}{j-1/3,n}\cup (\oU{jn}{0}\cap
\ol{T})$ in $M$, and to the neighborhood
\begin{equation*}
O_{jn}:=\left(\bigcap_{\alpha\in\N}\left\lfloor
\cU{jn\alpha}{jn},W_{jn}^{\alpha}
\right\rfloor\right)\cdot k_{jn}
\end{equation*}
of both $k_{jn}$ and $k'_{jn}$, where $k'_{jn}\in O_{jn}$ follows from
firstly $\Phi_{jn}(x)\in W_{jn}^{\alpha}$ and secondly
$k'_{jn}(x)=\Phi_{jn}(x)\cdot k_{jn}(x)\in W_{jn}^{\alpha}\cdot k_{jn}(x)$ for
all $x\in\cU{jn\alpha}{jn}$. Note that $O_{jn}$ is really open,
because Remark \ref{rem:locallyFiniteCoversByCompactSets} asserts that
just finitely many of the sets $\cU{jn\alpha}{jn}$ for fixed
$\alpha\in\N$ are non-empty and may influence the intersection.

By the choice of $U$, the result $\sm{k}_{jn}$ coincides with both
$k'_{jn}$ and $\sm{k}'_{jn}$ on $\bigcup_{i<j}\cU{ijn}{jn}$, so it
satisfies the cocycle conditions \eqref{eqn:cocycle}. It also
satisfies \eqref{eqn:quotient} by the choice of $O_{jn}$ and,
furthermore, \eqref{eqn:relative1}, because $\sm{k}_{jn}$ coincides
with $k'_{jn}$ on $\oU{jn}{0}\cap \ol{T}$ by the choice of $U$ and
$k'_{jn}$ coincides with $k_{jn}$ on $\oU{jn}{0}\cap T$.
\end{itemize}

This concludes the construction of the smooth principal $K$-bundle
$\cP_{\sm{k}}$. Since 
\[
\sm{k}_{jn}\big|_{\oU{jn}{\infty}\cap
T} = k_{jn}\big|_{\oU{jn}{\infty}\cap T},
\]
we may view $\left.\cP_{\sm{k}}\right|_{\eta^{-1}(T)}$ as the subset
$\left.\cP_{k}\right|_{\eta^{-1}(T)}$ of $\cP_{k}$. We use the same
covers of $M$ and identity neighborhoods in $K$ for the construction
of continuous functions $f_{i}:\cU{i}{0}\to K$ such that
\begin{enumerate}
\renewcommand{\labelenumi}{(\theenumi)}
\renewcommand{\theenumi}{\alph{enumi}}
\addtocounter{enumi}{3}
\item $f_{n}=\sm{k}_{nj}\cdot f_{j}\cdot k_{jn}$ pointwise on
$\cU{jn}{nn}$ for all $j<n$ in $\N$,\label{eqn:compatible}
\item $f_{n}\big(\cU{n\alpha}{0}\big)\subseteq W_{n}^{\alpha}$ for
all $\alpha,n\in\N$ and \label{eqn:small}
\item $f_{n}\big(\cU{n}{0}\cap T\big)=\{e\}$ for all
$n\in\N$.\label{eqn:relative2}
\end{enumerate}
Then Remark \ref{rem:bundleEquivalences} tells us that the restriction
of the maps $f_{i}$ to the sets $\oU{i}{\infty}$ of the open cover is
the local description of a bundle equivalence
$\Omega:P_{k}\to P_{\sm{k}}$. Indeed, all
the sets $\cU{jn}{nn}$ of condition \eqref{eqn:compatible} contain the
corresponding sets $\oU{jn}{\infty}$ of the open cover. Furthermore,
condition \eqref{eqn:relative2} implies that the restriction of
$\Omega$ to $\eta ^{-1}(T)$ is the identity.

We start with the constant function $f_{1}\equiv e$, which clearly
satisfies conditions \eqref{eqn:small} and \eqref{eqn:relative2}. Then
we construct $f_{n}$ for $n>1$ inductively as follows:
\begin{itemize}
\item To satisfy condition \eqref{eqn:compatible}, we start with
\begin{equation*}
f'_{n}:\bigcup_{j<n}\cU{jn}{jn}\to K,\quad
f'_{n}(x)=\sm{k}_{nj}(x)\cdot f_{j}(x)\cdot k_{jn}(x)\text{ for }
x\in\cU{jn}{jn}.
\end{equation*}
This continuous function is well-defined, because the conditions
\eqref{eqn:compatible} for $f_{j}$ on
$\cU{j'jn}{jn}\subseteq\cU{j'j}{jj}$ and \eqref{eqn:cocycle} for
$j'<j<n$ on $\cU{j'jn}{jn}$ guarantee that
\begin{equation*}
\sm{k}_{nj}(x)\cdot f_{j}(x)\cdot k_{jn}(x)
=\sm{k}_{nj}(x)\cdot\sm{k}_{jj'}(x)\cdot f_{j'}(x)
\cdot k_{j'j}(x)\cdot k_{jn}(x)
=\sm{k}_{nj'}(x)\cdot f_{j'}(x)\cdot k_{j'n}(x)
\end{equation*}
holds for all $x\in\cU{j'jn}{jn}$. In addition, $f'_{n}$ is the
identity on $\bigcup_{j<n}\cU{jn}{jn}\cap T$ by conditions
\eqref{eqn:relative1} and \eqref{eqn:relative2}.
\item To apply Lemma \ref{lem:fadingOut}, we need to know something
about the values of $f'_{n}$. For arbitrary $\alpha\in\N$ and
$x\in\cU{jn\alpha}{jn}$, conditions \eqref{eqn:quotient},
\eqref{eqn:small}, and \eqref{eqn:newf} of Lemma
\ref{lem:SuitableIdentityNeighborhoods} imply
\begin{align*}
f'_{n}(x)
&=\sm{k}_{nj}(x)\cdot f_{j}(x)\cdot k_{jn}(x) 
=k_{nj}(x)\cdot\big(\sm{k}_{jn}(x)\cdot k_{nj}(x)\big)^{-1}
\cdot f_{j}(x)\cdot k_{jn}(x) \\
&\in k_{nj}(x)\cdot\big(W_{jn}^{\alpha})^{-1}\cdot W_{j}^{\alpha}\cdot
k_{jn}(x) \subseteq W_{n}^{\alpha},
\end{align*}
so that the values of $f'_{n}$ are, altogether, contained in the
identity neighborhood $W$ of $K$. If we apply Lemma
\ref{lem:fadingOut} to $M:=\cU{n}{0}$, to $f'_{n}$ on
$A:=\bigcup_{j<n}\cU{jn}{jn}$ and to the smaller set
$B:=\bigcup_{j<n}\cU{jn}{nn}$, then we get a continuous function
$f_{n}:\cU{n}{0}\to W$, which satisfies the conditions \eqref{eqn:compatible},
\eqref{eqn:small} and \eqref{eqn:relative2}.
\end{itemize}
This concludes the construction of the bundle equivalence.
\end{prf}

\begin{theorem}[Smoothing Continuous Bundle Equivalences]
\label{thm:smoothingContinuousBundleEquivalences} Let $\cP$ and $\cP'$
be smooth principal $K$-bundles over the fi\-nite-\-di\-men\-sion\-al
connected paracompact manifold with corners $M$ and let $\Omega\from
P\to P'$ be a continuous bundle equivalence. If $C\se M$ is closed and
$\Omega$ is smooth on an open neighborhood of $\eta^{-1}(C)$, then
there exists an open neighborhood $T$ of $C$ and a smooth bundle
equivalence $\sm{\Omega}:P\to P'$ with
$\Omega\big|_{\eta^{-1}(T)}=\sm{\Omega}\big|_{\eta^{-1}(T)}$.
\end{theorem}

\begin{prf}
Let $(U_{j})_{j\in J}$ be an open cover of $M$ that is locally trivial
for both bundles $\cP$ and $\cP'$.  Proposition
\ref{prop:NestedCovers} yields locally finite open covers
$\big(\oU{i}{\lambda}\big)_{i\in\N}$ of $M$ for every
$\lambda\in\{0,\infty\}\cup\left(1+\frac{1}{3}\N\right)$ such that
the closures $\cU{i}{\lambda}$ are compact manifolds with corners and 
\begin{equation*}
\cU{i}{\infty }\se\oU{i}{j+1}\se \cU{i}{j+1}\se
\oU{i}{j+2/3}\se\cU{i}{j+2/3}\se\oU{i}{j+1/3}\se\cU{i}{j+1/3}\se
\oU{i}{j}\se \oU{i}{0}\se \cU{i}{0}\se U_{i}
\end{equation*}
holds for all $i,j\in \N$, where $U_{i}$ denotes a suitable set of the
cover $(U_{j})_{j\in J}$ for every $i\in\N$.  According to Remarks
\ref{rem:continuousBundle} and \ref{rem:smoothStructureOnBundle}, we
may then describe the smooth bundles $\cP$ and $\cP'$ by smooth
transition functions $k=(k_{ij})_{i,j\in\N}$ and
$k'=(k'_{ij})_{i,j\in\N}$ on the open cover $(U_{i})_{i\in\N}$,
equivalently, by their restrictions to any open cover
$\big(\oU{i}{\lambda}\big)_{i\in\N}$ from above. In these local
descriptions of the bundles, the bundle equivalence $\Omega$ can, as in
Remark \ref{rem:bundleEquivalences}, be seen as a family
$f_{i}:U_{i}\to K$ of continuous maps for $i\in\N$ that satisfy
\begin{equation}\label{eqn:compatibilityConditionOfContinuousEquivalence}
f_{j}(x)=\kp{ji}(x)\cdot f_{i}(x)\cdot \kq{ij}(x)
\fa i,j\in\N\text{ and }x\in U_{ij}.
\end{equation}
Let $S\se M$ be an open neighborhood of $C$ such that
$\left.\Omega\right|_{S}$ is smooth. This means that the restriction
of each $f_{i}$ to $U_{i}\cap S$ is smooth. In addition, there exists
an open neighborhood $T$ of $C$ with $\ol{T}\se S$ since $M$ is normal
(see Remark \ref{rem:topologicalPropertiesOfParacompactSpaces}).

We shall inductively construct smooth maps $\wt{f}_{i}\from
\cU{i}{0}\to K$ such that
\begin{enumerate}
\renewcommand{\labelenumi}{(\theenumi)}
\renewcommand{\theenumi}{\alph{enumi}}
\item $\sm{f}_{j}=k'_{ji}\cdot\sm{f}_{i}\cdot k_{ij}$ pointwise
on $\cU{ij}{j}$ for all $i<j$ in $\N$,
\label{eqn:compatibilityConditionInSmoothingEquivalences}
\item $\big(\sm{f}_{i}\cdot f^{-1}_{i}\big)\big(\cU{i\alpha}{i}\big)\subseteq
W_{i}^{\alpha}$ for all $i,\alpha\in\N$ and
\label{eqn:quotientConditionInSmoothingEquivalences}
\item $\left.\sm{f}_{i}\right|_{\cU{i}{0}\cap
T}=\left.{f}_{i}\right|_{\cU{i}{0}\cap T}$ for all $i\in\N$
\label{eqn:relative3}
\end{enumerate}
are satisfied at each step, where the $W_{i}^{\alpha}$ are
$\varphi$-convex identity neighborhoods provided by Lemma
\ref{lem:SuitableIdentityNeighborhoods} that we apply to the countable
compact cover $\big(\cU{i}{0}\big)_{i\in\N}$, to the transition
functions $k'_{ij}$, to a convex centered chart
$\varphi:W\to\varphi(W)$ of $K$ and to the sequence of unit
neighborhoods which is constantly $W$ (we do not need the
$W_{ij}^{\alpha}$ in this proof). These maps $\sm{f}_{i}$ describe a
smooth bundle equivalence between $\cP$ and $\cP'$ when restricted to
the open cover $\big(\oU{i}{\infty}\big)_{i\in\N}$, because
$\eqref{eqn:compatibilityConditionInSmoothingEquivalences}$ asserts
that $\sm{f}_{j}=k'_{ji}\cdot\sm{f}_{i}'\cdot k_{ij}$ is satisfied on
$\oU{ij}{\infty}$ for all $i<j$, in particular.

To construct the smooth function $\wt{f}_{1}:\cU{1}{0}\to K$, we apply
Proposition \ref{prop:smoothing} to the continuous map $f:=f_{1}$ on
$M:=\cU{1}{0}$, the closed set $A:=\cU{1}{0}\cap \ol{T}$, the open set
$U:=\cU{1}{0}\backslash \ol{T}$ and to the open neighborhood
\begin{equation*}
O_{1}:=\bigcap_{\alpha \in \N}
\left\lfloor\cU{1\alpha}{0},W_{1}^{\alpha}\right\rfloor\cdot f_{1}
\end{equation*}
of $f_{1}$, which is indeed open, since only finitely many
$\cU{1\alpha}{0}$ are non-empty by Remark
\ref{rem:locallyFiniteCoversByCompactSets}.  By construction,
$\sm{f}_{1}$ satisfies
\eqref{eqn:quotientConditionInSmoothingEquivalences} and
\eqref{eqn:relative3}.  To construct the smooth function
$\wt{f}_{j}:\cU{j}{0}\to K$ inductively for $j>1$, we need the usual
three steps:
\begin{itemize}

\item In order to satisfy
\eqref{eqn:quotientConditionInSmoothingEquivalences} in the end, we
define a continuous map
\begin{equation*}
\sm{f}'_{j}\from\bigcup_{i<j}\cU{ij}{j-1}\to K,\quad
\wt{f}_{j}'(x):=\kp{ji}(x)\cdot \wt{f}_{i}(x)\cdot \kq{ij}(x)
\text{ for }x\in \cU{ij}{j-1}.
\end{equation*}
If $x$ is an element of both $\cU{ij}{j-1}$ and $\cU{i'j}{j-1}$ for
$i'<i<j$, condition
\eqref{eqn:compatibilityConditionInSmoothingEquivalences} for $j-1$
and the cocycle conditions of both $k$ and $k'$ assert that the
so-defined values for $\sm{f}_{j}'(x)$ agree. Furthermore, the
compatibility condition together with condition \eqref{eqn:relative3}
ensure that $\sm{f}'_{j}$ coincides with $f_{j}$ on
$\bigcup_{i<j}\cU{ij}{j-1}\cap T$.

\item This definition of $\sm{f}_{j}'$, along with
\eqref{eqn:compatibilityConditionOfContinuousEquivalence} and property
\eqref{eqn:newf} in Lemma \ref{lem:SuitableIdentityNeighborhoods}
assert that
\begin{equation*}
\varphi_{j}(x):= \wt{f}_{j}'(x)\cdot f_{j}(x)^{-1}= \kp{ji}(x)\cdot
\wt{f}_{i}(x)\cdot \kq{ij}(x)\cdot f_{j}(x)^{-1}= \kp{ji}(x)\cdot
\underbrace{\wt{f}_{i}(x)\cdot f_{i}(x)^{-1}}_{\in W_{i}^{\alpha
}}\cdot \kp{ij}(x) \in W_{j}^{\alpha}
\end{equation*}
holds for all $x\in \cU{ij\alpha}{j-1}$, $i<j$ and $\alpha$ in
$\N$. So we may apply Lemma \ref{lem:fadingOut} to
$A\nobreak:=\nobreak\bigcup_{i<j}\cU{ij}{j-1}$ and
$B:=\bigcup_{i<j}\cU{ij}{j-2/3}$ to fade out $\varphi_{j}$ to a
continuous map $\Phi_{j}$ on $M:=\cU{j}{0}$.  Then $\Phi_{j}$
coincides with $\varphi_{i}$ on $B$ and maps $\cU{j\alpha}{j}$ into
$W_{j}^{\alpha}$. Since $\sm{f}'_{j}$ coincides with $f_{j}$ on $A\cap
T$, $\varphi_{j}$ and, consequently, $\Phi_{j}$ is the identity on
$(A\cap T)\cup (\cU{j}{0}\backslash A)$.  Furthermore, $\varphi_{j}$
is smooth on $A^{0}\cap S$, because so are all its constituents, and
on a neighborhood of $\cU{j}{0}\setminus A^{0}$, because it is the
identity there. By the last conclusion of Lemma \ref{lem:fadingOut},
$\Phi_{j}$ is smooth on $\cU{j}{0}\cap S$.

\item Accordingly, $\Phi_{j}\cdot f_{j}$ is an element of the open
(due to Remark \ref{rem:locallyFiniteCoversByCompactSets})
neighborhood
\begin{equation*}
O_{j}:=\bigcap_{\alpha\in\N}\left\lfloor\cU{j\alpha}{j},W_{j}^{\alpha
}\right\rfloor\cdot f_{j}
\end{equation*}
of $f_{j}$ and is smooth on $\bigcup_{i<j}\oU{ij}{j-2/3}$ and on
$\cU{j}{0}\cap S$. If we apply Proposition \ref{prop:smoothing} to
$M:=A:=\cU{j}{0}$, $U:=M\setminus\left(\ol{T}\cup
\bigcup_{i<j}\cU{ij}{j-1/3}\right)$, $O_{j}$, and to $f:=\Phi_{j}\cdot
f_{j}$, then we obtain a smooth map $\wt{f}_{j}:\cU{j}{0}\to K$.
\end{itemize}
The map $\sm{f}_{j}$ satisfies
\eqref{eqn:compatibilityConditionInSmoothingEquivalences}, because so
does $\sm{f}_{j}'$, with which it coincides on
$\bigcup_{i<j}\cU{ij}{j}$. Moreover,
\eqref{eqn:quotientConditionInSmoothingEquivalences} is satisfied due
to the choice of $O_{j}$. So is \eqref{eqn:relative3}, because
$\Phi_{j}$ is the identity on $\cU{j}{0}\cap T$ and $\Phi_{j}\cdot
f_{j}$ remains unchanged on $M\setminus U \supseteq \cU{j}{0}\cap T$ in the last step. This
concludes the construction.
\end{prf}

\begin{corollary}[Equivalences of Smooth and Continuous Bundles]
\label{cor:mainResult}
Let $K$ be a Lie group modeled on a locally convex space and $M$ be a
fi\-nite-\-di\-men\-sion\-al connected manifold with corners. Then
each continuous principal $K$-bundle over $M$ is continuously
equivalent to a smooth principal bundle. Moreover, two smooth
principal $K$-bundles over $M$ are smoothly equivalent if and only if
they are continuously equivalent.
\end{corollary}

\begin{prf}
The first statement is Theorem
\ref{thm:smoothingContinuousPrincipalBundles}, applied to
$C=\emptyset$. In the same way, the second assertion follows from
Theorem \ref{thm:smoothingContinuousBundleEquivalences}.
\end{prf}

\section{Related Topics}

In this section, we explain the relations of the results of the
preceding section to the problem of extending bundles from
submanifolds, to non-abelian \v{C}ech cohomology and to twisted
$K$-theory. One encounters the first situation, e.g., in
the construction of topological field theories. While non-abelian
\v{C}ech cohomology is only an equivalent sheaf-theoretic framework
for the problem, we show in the end how applications arise in twisted
$K$-theory.

\newcommand{\im}{\ensuremath{\mathop{\text{im}}}}

\begin{remark}[Abelian \v{C}ech Cohomology] Let $M$ be a paracompact
topological space with an open cover $\mathcal{U}=(U_{i})_{i\in I}$
and $A$ be an abelian topological group. Then for $n\geq 0$, an
\emph{$n$-cochain} $f$ is a collection of continuous functions
$f_{i_{1}\ldots i_{n+1}}:U_{i_{1}\ldots i_{n+1}}\to A$, and we denote
the set of \mbox{$n$-cochains} by $C^{n}(\mathcal{U},A)$ and set it to
$\{0\}$ if $n< 0$. We then define the boundary operator
\[
\delta_{n}\from C^{n}(\mathcal{U},A)\to C^{n+1}(\mathcal{U},A),\;\;
\delta (f)_{i_{0}i_{1}\ldots i_{n+1}}
=\sum_{k=0}^{n}(-1)^{k}f_{i_{0}\ldots\wh{i_{k}}\ldots i_{n+1}},
\]
where $\wh{i_{k}}$ means that we omit the index $i_{k}$. Then
$\delta_{n+1}\circ \delta_{n}=0$, and we define
\begin{equation}\label{eqn:definitionOfAbelianCechCohomologyGroups}
\check{H}_{c}^{n}(\mathcal{U},A):=\ker(\delta_{n})/\im(\delta_{n-1})
\quad\text{ and }\quad 
\check{H}_{c}^{n}(M,A):=\lim_{\to}\;\check{H}_{c}^{n}(\mathcal{U},A).
\end{equation}
The group $\check{H}^{1}(M,A)$ is the \textit{$n$-th continuous \v{C}ech
cohomology}. If, in addition, $M$ is a smooth manifold with or without
corners and $A$ is a Lie group, then the same construction with
smooth instead of continuous functions leads to the corresponding
\textit{$n$-th smooth \v{C}ech cohomology}.
\end{remark}

\begin{theorem}[Isomorphism for Abelian \v{C}ech Cohomology]
\label{thm:smoothingHigherCocycles} Let $M$ be a
fi\-nite-\-di\-men\-sion\-al connected manifold and $A$ be an abelian
locally convex Lie group, then for each $n\in \N$, the canonical
map $\iota \from \check{H}^{n}_{s}(M,A)\to \check{H}^{n}_{c}(M,A)$
defines an isomorphism of abelian groups.
\end{theorem}

\begin{prf}
It clearly suffices to show that $\iota$ is bijective, so take some
$[(f_{i})_{i\in \N^{n}}]$, defining an element of $H^{n}_{c}(M,A)$.
Then $[(f_{i})_{i\in \N^{n}}]$ is in some $H^{n}_{c}(\mathcal{U},A)$
and as before, we may assume that $\mathcal{U}$ is a countable cover
of $M$. Choosing a bijection $\lambda \from \N^{n}\to \N^{+}$ induces
a total order on $\N^{n}$. From Proposition \ref{prop:NestedCovers},
we obtain corresponding open covers $\mathcal{U}^{[i]}$ for each $i\in
\N^{n}\cup\{0,\infty\}$ such that
$\cU{i}{\infty}\se\oU{i}{k}\se\cU{i}{k}\se \oU{i}{k'}\se \cU{i}{k'}\se
\oU{i}{0}$ if $k'<k$. In addition, we choose an arbitrary
$\varphi$-convex neighbourhood $W$ in $A$.  We set the stage for the
induction by defining $\wt{f}_{i}$ for the least $n-1$ elements
$i_{1},\dotsc ,i_{n-1}$ of $\N^{n}$ to be smooth functions
$\wt{f}_{i}\in C^{\infty}(U_{i},A)$ with $(\wt{f_{i}}\cdot
f_{i}^{-1})(U_{i})\se W$ (cf.\ Proposition \ref{prop:smoothing}).

Defining $\wt{f}_{k}$ inductively on $\cup_{i<k}(U_{i}\cap U_{k})$ by
the cocycle condition, fading it out appropriately to $U_{k}$ and
smoothing it out we obtain $\wt{f_{k}}\in C^{\infty}(U_{k},A)$,
satisfying $(\wt{f}_{k}\cdot f_{k}^{-1})(U_{k})\se W$ and satisfying
the cocycle condition on $(\oU{i}{k})_{i\in \N}$. In the end, this
yields a cocycle $[\wt{f}_{i}]$ in
$\check{H}^{n}_{s}(\mathcal{U}^{[\infty]},A)$, which is equivalent to
$[(\left.f_{i}\right|_{\oU{i}{\infty}})_{i\in \N^{n}}]$. This yields
the surjectivity and the injectivity may be achieved analogously.
\end{prf}

\begin{remark}[On the previous proof] Note we have no occurrence of
Lemma \ref{lem:SuitableIdentityNeighborhoods} in the abelian case, as
the adjoint action, occurring implicitly in the cocycle condition, is
trivial in this case. Thus we can achieve that $\wt{f}_{k}\cdot
f_{k}^{-1}$ always has values in one fixed identity neighbourhood,
simplifying the proof significantly.
\end{remark}

\begin{remark}[Non-Abelian \v{C}ech Cohomology] (cf.\ \cite[Sect.\
12]{dedecker54} and \cite[3.2.3]{gelfandManin99}) If $n=0,1$, then we
can perform a similar construction as in the previous remark in the
case of a not necessarily abelian group $K$. The definition of an
\emph{$n$-cochain} is the same as in the commutative case, but we run
into problems when writing down the boundary operator $\delta$.
However, we may define $\delta_{0}(f)_{ij}=f_{i}\cdot f_{j}^{-1}$,
$\delta_{1}(k)_{ijl}=k_{ij}\cdot k_{jl}\cdot k_{li}$ and call the
elements of $\delta_{1}^{-1}(\{e\})$ \emph{$2$-cocycles} (or
\emph{cocycles}, for short).

The way to circumvent difficulties for $n=1$ is the observation
that even in the non-abelian case, $C^{1}_{c}(\mathcal{U},K)$ acts on
cocycles by $(f_{i},k_{ij})\mapsto f_{i}\cdot k_{ij}\cdot
f_{j}^{-1}$. Thus we define two cocycles $k_{ij}$ and $k_{ij}'$ to be
equivalent if $k_{ij}'=f_{i}\cdot k_{ij}\cdot f_{j}^{-1}$ on $U_{ij}$
for some $f_{i}\in C^{1}(\mathcal{U},K)$, and by
$\check{H}_{c}^{1}(\mathcal{U},K)$ the equivalence classes (or the orbit space)
of this action. Then $\check{H}_{c}^{1}(\mathcal{U},K)$ is not a group, but we
may nevertheless take the direct limit
\[
\check{H}_{c}^{1}(M,K):=\lim_{\to}\; \check{H}_{c}^{1}(\mathcal{U},K)
\]
of sets and define it to be the \emph{$1^{\text{st}}$ (non-abelian)
continuous \v{C}ech cohomology} of $M$ with coefficients in $K$. A
representing space of $\check{H}_{c}^{1}(M,K)$ would then be the set of
equivalence classes of continuous principal $K$-bundles over $M$.

Again, if $M$ is a smooth manifold with corners and $K$ is a
Lie group, we can adopt this construction to define the
\emph{$1^{\text{st}}$ (non-abelian) smooth \v{C}ech cohomology}
$\check{H}^{1}_{s}(M,K)$.
\end{remark}

\begin{theorem}[Isomorphism for Non-Abelian \v{C}ech Cohomology]
\label{thm:isomorphismForNonAbelianCechCohomology}
If $M$ is a fi\-nite-di\-men\-sion\-al connected manifold with corners and
$K$ is a Lie group modeled on a locally convex space, then the
canonical map
$\iota :\check{H}_{s}^{1}(M,K)\to \check{H}^{1}_{c}(M,K)$
is a bijection.
\end{theorem}

\begin{prf}
We identify smooth and continuous principal bundles with \v{C}ech
$1$-cocycles and smooth and continuous bundle equivalences with
\v{C}ech $0$-cochains as in Remark \ref{rem:bundleEquivalences}.  For
each open cover $\mathcal{U}$ of $M$, we have the canonical map
$\check{H}^{1}_{s}(\mathcal{U},K)\to
\check{H}^{1}_{c}(\mathcal{U},K)$. Now each cocycle $k_{ij}\from
U_{ij}\to K$ defines a principal bundle $\cP$ with locally trivial
cover $\mathcal{U}$. We may assume by Corollary \ref{cor:mainResult}
that $\cP$ is continuously equivalent to a smooth principal bundle
$\wt{\cP}$, and thus that $\mathcal{U}$ is also a locally trivial
covering for $\wt{\cP}$. This shows that the map is surjective, and
the injectivity follows from Corollary \ref{cor:mainResult} in the
same way. Accordingly, the map induced on the direct limit is a
bijection.
\end{prf}

\renewcommand{\bH}{\mathcal{H}}

\begin{remark}[The Projective Unitary Group] Let $\bH$ be a separable
infi\-nite-\-di\-men\-sion\-al Hilbert space and denote by
$\operatorname{U}(\bH)$ the group of unitary operators. If we equip
$\operatorname{U}(\bH)$ with the norm topology, then the exponential
series, restricted to skew-self-adjoint operators
$L(\operatorname{U}(\bH))$, induces a Banach--Lie group structure on
$\operatorname{U}(\bH)$ (cf.\ \cite[Ex.\ 1.1]{milnor84}). Then
$\operatorname{U}(1)\cong Z(U(\bH ))$,
and it can also be shown that
$\operatorname{PU}(\bH):=\operatorname{U}(\bH)/\operatorname{U}(1)$ is
a Lie group modeled on $L(\operatorname{U}(\bH))/i\R$.
\end{remark}

\begin{remark}[Eilenberg--MacLane Spaces] If $X$ is a topological
space with trivial $n$-th homotopy group $\pi_{n}(X)$ for all but
one $n\in \N$, then it is called an \textit{Eilenberg--MacLane space}
$K(n,\pi_{n}(X))$.  Since $\operatorname{U}(1)$ is a $K(1,\Z)$, the
long exact homotopy sequence \cite[Th.\ VII.6.7]{bredon93} shows that
$\operatorname{PU}(\bH)$ is a $K(2,\Z)$, since $\operatorname{U}(\bH)$
is contractible \cite[Th.\ 3]{kuiper65}. By the same argument, the
classifying space $B\,\operatorname{PU}(\bH)$
is a $K(3,\Z)$, since its total space
$E\,\operatorname{PU}(\bH)$ is contractible. Thus
\[
\check{H}^{3}(M,\Z)\cong [M,B\,\operatorname{PU}(\bH)]\cong
\check{H}_{c}^{1}(M,\operatorname{PU}(\bH))
\]
by \cite[Cor.\ VII.13.16]{bredon93}. The representing class $[\cP]$ in
$\check{H}^{3}(M,\Z)$ is called the \textit{Dixmier-Douady class} of
$\cP$ (cf.\ \cite{careyCrowleyDiarmuidMurray98},
\cite{dixmierDouady63}).  It describes the obstruction of $\cP$ to be
the projectivization of an (automatically trivial) principal
$\operatorname{U}(\bH)$-bundle.
\end{remark}

\begin{corollary}[Smoothing
$\boldsymbol{\operatorname{PU}(\bH)}$-bundles] If $M$ is a
fi\-nite-\-di\-men\-sion\-al connected manifold with corners, then
\[
\check{H}^{3}(M,\Z)\cong \check{H}^{1}_{c}(M,\operatorname{PU}(\bH ))
\cong \check{H}^{1}_{s}(M,\operatorname{PU}(\bH )).
\]
\end{corollary}

\newcommand{\FrH}{\ensuremath{\mathop{\text{Fred}(\bH)}}}
\newcommand{\Ad}{\ensuremath{\mathop{\text{Ad}}}}

\begin{remark}[Twisted $\mathbf{K}$-Theory] (cf.\ \cite[Sect.\
2]{rosenberg89}, \cite{mathaiMurray02}) The Dixmier-Douady class of a
principal $\operatorname{PU}(\bH)$-bundle over $M$ induces a twisting
of the $K$-theory of $M$ in the following manner. For any paracompact
space, the $K$-theory $K^{0}(M)$ is defined to be the Grothendieck
group of the monoid of equivalence classes of
fi\-nite-\-di\-men\-sion\-al complex vector bundles over $X$, where
addition and multiplication is defined by taking direct sums and
tensor products of vector bundles \cite{husemoller94}. Furthermore,
the space of Fredholm operators $\FrH$ is a representing space for
$K$-theory, i.e. $K^{0}(M)\cong [M,\FrH]$, where $[\cdot ,\cdot ]$
denotes homotopy classes of continuous maps.  Since
$\operatorname{PU}(\bH)$ acts (continuously) on $\FrH$ by conjugation,
we can form the associated vector bundle $\cP_{\FrH}:=\FrH
\times_{\operatorname{PU}(\bH)}\cP $.  Then the homotopy classes of
sections $[M,P_{\FrH}]$ (or equivalently, the equivariant homotopy
classes of equivariant maps
$[P_{\FrH},\FrH]^{\operatorname{PU}(\bH)}$) define the \textit{twisted
$K$-theory} $K_{\cP}(M)$.  Now Corollary \ref{cor:mainResult} implies
that we may assume $\cP$ to be smooth. Since the action of
$\operatorname{PU}(\cH)$ on $\FrH$ is locally given by conjugation, it
is smooth (recall that $\FrH$ is an open subset of $B(\cH)$,
giving it a natural manifold structure), whence is $\cP_{\FrH}$. Due
to Lemma \ref{lem:smoothingHomotopies}, we may, in the computation of
$K_{\cP}(M)$, restrict our attention to smooth sections and smooth
homotopies.
\end{remark}

\section*{Acknowledgements} The authors of this paper are grateful to
Mathai Varghese and Karl-Hermann Neeb for pointing out the problem and
to Peter Teichner for pointing out the need of Theorem
\ref{thm:smoothingHigherCocycles}.

\def\polhk#1{\setbox0=\hbox{#1}{\ooalign{\hidewidth
  \lower1.5ex\hbox{`}\hidewidth\crcr\unhbox0}}} \def\cprime{$'$}

\vskip\baselineskip
\vskip\baselineskip
\vskip\baselineskip
\large
\noindent 
Christoph M\"uller\\
Christoph Wockel\\
Fachbereich Mathematik\\
Technische Universit\"at Darmstadt\\
Schlossgartenstra\ss e 7\\
D-64289 Darmstadt\\
Germany\\[\baselineskip]
\normalsize
\texttt{cmueller@mathematik.tu-darmstadt.de}\\
\texttt{christoph@wockel.eu}
\end{document}